\documentclass[11pt]{amsart}

\ifx\textheight\undefined
\else
\advance\voffset by1.2cm
\advance\voffset by-2.3cm

\usepackage{amssymb}
\usepackage{amscd}
\usepackage{amsthm}

\newtheorem{Theorem}{Theorem}[section]
\newtheorem{Lemma}[Theorem]{Lemma}
\newtheorem{Corollary}[Theorem]{Corollary}
\newtheorem{Proposition}[Theorem]{Proposition}

\newtheorem{Question}[Theorem]{Question}

\newcommand{\R}{\mathbb{R}}

\newcommand{\tr}{\operatorname{tr}}

\newcommand{\Irr}{\operatorname{Irr}}

\newcommand{\kernel}{\operatorname{Ker}}

\newcommand{\SL}{{\rm SL}}
\newcommand{\Sym}{{\rm S}}
\newcommand{\Alt}{{\rm A}}

\begin{document}

\title{A local conjecture on Brauer character degrees of finite
groups}

\author{Thorsten Holm, Wolfgang Willems}

\address{Thorsten Holm\newline
Department of Pure Mathematics,
University of Leeds \newline
Leeds LS2 9JT, United Kingdom}
\email{tholm@maths.leeds.ac.uk}

\address{Wolfgang Willems\newline
Otto-von-Guericke-Universit\"at,
Institut f\"ur Algebra und Geometrie\newline
Postfach 4120,
39016 Magdeburg,
Germany}
\email{wolfgang.willems@mathematik.uni-magdeburg.de}


\begin{abstract}
Recently, a new conjecture on the degrees of the irreducible 
Brauer characters of a finite group was presented in
\cite{Willems}.
In this paper we propose a 'local' version of this conjecture for
blocks $B$ of finite groups, giving 
a lower bound for $\sum \varphi(1)^2$ where the sum runs through the set of irreducible
Brauer characters of $B$
 in terms of invariants of $B$.
A slight reformulation leads to 
interesting open questions about traces of Cartan matrices
of blocks. 

We show that the local conjecture is true for blocks with
one simple module, blocks of
$p$-solvable groups and blocks with cyclic defect groups.
It also holds for many further examples of blocks of sporadic groups,
symmetric groups or groups of Lie type. 
Finally we prove that the  conjecture is true for 
blocks of tame representation type.
\end{abstract}

\maketitle


\section*{Introduction}


Representation theory of finite groups is an area which is
characterized by an enormous number of important long-standing
open problems. Many of them deal with ordinary or 
modular characters, so
character theory still is one of the most important and fundamental 
parts of representation theory. The key aim 
in this area is to find intimate relations between 
the numerical invariants given by characters and
the structure of the group. Actually,
R. Brauer's main motivation for introducing
modular characters was to study the structure of finite simple groups. 

In this paper we are going to study some of the most
important numerical invariants in modular representation
theory, namely degrees of irreducible Brauer
characters. The main aim is to propose a conjecture relating
the degrees of the irreducible characters belonging to a block 
to some well-known invariants of the block. This can be considered as
a block version of a conjecture recently presented in 
\cite{Willems} for finite groups. In order to state the conjecture 
we need to introduce some notation.

Let $G$ be a finite group and let $p$ be a prime dividing the 
order of $G$. Moreover, 
let $B$ be a $p$-block of $G$ with 
defect group $D$ where the
underlying field of characteristic $p$ is assumed to be algebraically 
closed.
Furthermore let $\Irr_p(B) = \{ \varphi_1, \ldots, 
\varphi_{l(B)} \}$ denote the set of irreducible
Brauer characters belonging to $B$, and $\Irr_p(G)$ the set of all 
irreducible Brauer
characters of $G$. 

In \cite{Willems} the second author proposed the following conjecture
relating the degrees of the irreducible Brauer characters
of a finite group $G$ to the $p'$-part of the order of $G$
(i.e. the largest divisor of the order of $G$ not divisible by 
$p$). 
\medskip

\noindent
{\bf Global Conjecture} \cite{Willems}.
{\em For any finite group $G$ we have
\begin{eqnarray} \label{nr1a}
|G|_{p'} \leq  \sum_{\varphi \in \Irr_p(G)} \varphi(1)^2
\end{eqnarray}
with equality if and only if the Sylow $p$-subgroup of $G$ is 
normal. }
\medskip

This conjecture is known to be true for $p$-solvable groups,
groups with a cyclic Sylow $p$-subgroup, groups of Lie type in the 
defining characteristic, and 
asymptotically in many cases in non-defining characteristic, 
and many more examples (see \cite{Willems} for details). 

The main aim of this paper is to develop a `local' version of this
conjecture, i.e. a version of this conjecture for blocks of 
finite groups. The crucial aspect is to find a good replacement for
the left-hand side, which only depends on the block and its 
invariants, but not on the group to which the block belongs. 

We propose the following local version,
relating the degrees of the irreducible Brauer characters to 
well-known invariants of the block. By $ \dim B$ we denote the dimension
of a block $B$ as vector space over the ground field.
\medskip

\noindent
{\bf Local Conjecture.} 
{\em Let $B$ be a $p$-block of a finite group, and let
$D$ be its defect group. Then we have
\begin{eqnarray} \label{nr1b} 
  \frac{\dim B}{l(B)|D|} \leq  \sum_{\varphi \in \Irr_p(B)} \varphi(1)^2
 \end{eqnarray}
with equality if and only if $l(B)=1.$} 
\medskip

When experimentally computing examples one observes that in most 
cases the right-hand side becomes far larger than the left-hand side,
and at first glance one is tempted to replace the right-hand 
side of (\ref{nr1b}) by 
the stronger bound
$$ \max_{\varphi \in \Irr_p(B)} \varphi(1)^2.$$
However, as Christine Bessenrodt pointed out to us, 
this does not hold true in general,
as the next example shows.\\

\noindent
{\bf Example} Let $G= \Sym_{10}$ be the symmetric group on $10$ letters 
and let $p=2$.
By \cite{Kerber}, the group $G$ has a $2$-block of defect $3$ 
with two irreducible
Brauer characters of degrees $128$ and $160$, respectively. 
The degrees of the irreducible
classical characters are $160, 160, 448, 288$ and $288$. Thus
$$ \frac{\dim B}{l(B)|D|} = 26112,$$
but
$$ \max_{\varphi \in \Irr_p(B)} \varphi(1)^2 = 25600.$$

\medskip

We would like to mention that there are other examples in 
characteristic $2$,
but we are not aware of any such example in odd characteristic.
\smallskip

Denote by $l(G)$ the 
number of $p'$-conjugacy classes of $G$ (i.e. conjugacy classes of
$G$ of elements whose order is not divisible by $p$). By a classical
result of Brauer, $l(G)$ equals the number of irreducible
Brauer characters of $G$. With this notation the following is obviously 
a weak version of the Global Conjecture.  
\medskip

\noindent
{\bf Weak Global Conjecture.} 
{\em For any finite group $G$ we have 
$$
  \frac{|G|_{p'}}{l(G)} \leq  \max_{\varphi \in \Irr_p(G)} \varphi(1)^2
$$
with equality if and only if $G=O_p(G) \times O_{p'}(G) $, where 
$O_{p'}(G)$ is abelian.
} 
\medskip

At this moment, in either form the conjectures
seem to be difficult to prove in general. 
Still, in this paper we are going 
to show that the Local Conjecture is true for various 
important classes of blocks of finite groups. 
The results can be 
summarized as follows.
\medskip

\noindent
{\bf Theorem 1.} {\em 
Let $B$ be a $p$-block of a finite group $G$ with defect group
$D$, and let $l(B)$ be the number of irreducible Brauer characters.
Then we have
\begin{eqnarray} \label{strong}
\frac{\dim B}{l(B)|D|} \leq  \max_{\varphi \in \Irr_p(B)} 
    \{\varphi(1)^2\}
\end{eqnarray}     
if one of the following holds:
\begin{enumerate}
\item[{(a)}] $l(B)=1$, with equality.
\item[{(b)}] $G$ is $p$-solvable.

\noindent
Equality holds in {\rm (\ref{strong})} if and only
if the following holds:
\begin{enumerate} 
\item[{(i)}] All irreducible Brauer characters of $B$ have the same
degree.
\item[{(ii)}] $\frac{|G|_p}{|D|}=\varphi(1)_p$ for all $\varphi\in 
\Irr_p(B)$. 
\end{enumerate}
\item[{(c)}] $D$ is cyclic.

\noindent
Equality holds in {\rm (\ref{strong})} if 
the corresponding Brauer tree is a star with exceptional vertex in the 
center. 
\end{enumerate}
In particular, in all these cases the Local Conjecture holds true.
}

\medskip

Moreover, 
in addition to these classes of blocks we checked many other
examples for which the Local Conjecture is true, giving further
evidence for its validity. 
Among these examples are many blocks 
of sporadic simple groups, blocks of (small) symmetric or 
alternating groups, and certain blocks of groups of Lie type. 
\smallskip

For proving the Local Conjecture we have the following interesting 
approach
via linear algebra. Namely, 
the following close relation to Cartan 
matrices will be proved in Section \ref{Sec-eigen}.
\medskip

\noindent
{\bf Proposition 2.} {\em Let $B$ be a $p$-block of a finite 
group, and  let $C$ be the 
 Cartan matrix of $B$ with trace ${\tr C}$.  
Then the following holds.\begin{itemize}
   \item[{(a)}] We have
$$\frac{\dim B}{\tr C} \leq \sum_{\varphi \in \Irr_p(B)} 
\varphi(1)^2, $$          
and equality holds if and only if $l(B)=1$.
    \item[{(b)}] If 
                 ${\tr C}\le l(B)|D|$,
                 then the Local Conjecture holds for $B$.
\end{itemize}
}

\medskip

This linear algebra approach seems to be helpful where the
Cartan invariants are known but not the degrees of the irreducible 
Brauer characters. In this way we can prove the
following result for blocks of tame representation type 
based on Erdmann's classification of their Cartan matrices. 
\medskip

\noindent
{\bf Theorem 3:} {\em The Local Conjecture holds for all blocks
of tame representation type (i.e. for all {\em 2}-blocks with defect group
dihedral, semidihedral or generalized quaternion).
}
\medskip

We are not aware of any example of a block such that
${\tr C}\le l(B)|D|$ does not hold. Actually, we strongly believe
that this inequality will be true in general. 
Note that the inequality is not at all obvious since there are
examples of blocks where the diagonal entries in the Cartan  
matrix are larger than $|D|$. The first example was given by 
Landrock \cite{Landrock}, disproving Problem 22 of Brauer's
famous list \cite{Brauer}.
Even more striking examples were found by Chastkofsky and Feit
\cite{Chast-Feit}; they showed that for the principal $2$-blocks
of the Suzuki groups
$\mbox{Suz}(2^m)$ the first diagonal Cartan entry is asymptotically
$|D|^{3/2}$. 

Further evidence for the validity of our global conjectures is provided by 
analogies between degrees of irreducible characters and lengths of
conjugacy classes. Such analogies have been observed by several authors 
(see for instance \cite{Huppert-Mainz}). 
The statements about conjugacy class lengths which correspond to 
our Global Conjecture and Weak Global Conjecture on character degrees
can actually be proved. We plan to address this issue in a separate 
note.


\section{Discussions around the conjectures} \label{Sec-why}


In this section we collect some of the evidence we have supporting
the Local Conjecture, and discuss various aspects around it.
\medskip

\noindent
{\bf 1.1 Blocks with one simple module.} 
Let $B$ be a block with one simple module, i.e. $l(B)=1$.
The validity of the Local Conjecture in this case is a classical 
result of R. Brauer, stating that 
in the case $l(B)=1$ equality holds in (\ref{nr1b}) 
(see \cite{Brauer76}, Theorems 2 and 4). This proves 
Theorem 1(a).

For an arbitrary block $B$, Brauer showed in \cite{Brauer76} that
the $p$-part of
$ \frac{\dim B}{|D|}$ divides the $p$-part of 
$\sum_{\varphi \in \Irr_p(B)} \varphi(1)^2$.
Thus equality in (\ref{nr1b}) forces $ p \nmid l(B)$.
But the `only if' part of the characterization of equality
remains open. 
\medskip

\noindent
{\bf 1.2} It might be tempting to try to improve the left-hand side
of the inequality in the Local Conjecture. We shall see below that for
some interesting classes of blocks one 
can actually show that  
\begin{eqnarray} \label{nr1e}
\frac{\dim B}{|D|} \le \sum_{\varphi \in \Irr_p(B)} \varphi(1)^2,
\end{eqnarray}
which obviously is stronger than the statement in the 
Local Conjecture for $B$. 
\smallskip

But the following examples show that the 
factor $l(B)$ in the Local Conjecture can in general not be dropped and 
not be replaced by any constant. 

Let $G= \Alt_5$ denote the alternating group on 
$5$ letters and
let $B$ be the principal $2$-block of $G$. Then $B$ contains 
three irreducible
Brauer characters of degrees $1,2$ and $2$, respectively.
The corresponding projective 
characters are of degree $12,8,8$. Thus we obtain
$$ \frac{ \dim B}{|D|} =\frac{44}{4} = 11 > 9 = 1+4+4 = 
\sum_{\varphi \in \Irr_2(B)} \varphi(1)^2.$$

Now we take $n$ copies of $G=\Alt_5$ and consider the 
principal $2$-block $B^{(n)}$
of the direct product $G^{(n)}= G \times \ldots \times G$. 
Clearly $B^{(n)}$ has defect group $D^{(n)} = D \times \ldots \times D$.
It is not hard to compute that
$$ \frac{\dim B^{(n)}}{|D^{(n)}|} = 11^n > 9^n = \sum_{\varphi \in 
\Irr_2(B^{(n)})} \varphi(1)^2.$$
Thus the factor $l(B)$ in the Local Conjecture can not be replaced by a 
constant. 
\smallskip

Similar examples are provided by the maximal non-principal
$3$-block of  $6\Alt_7$, and by the principal $2$-block of $\SL(2,5)$,
showing that the above phenomenon depends on neither
the prime 2 nor on $D$ being abelian. 
\medskip


\noindent
{\bf 1.3}
We would like to mention that equality in (\ref{nr1e}) does not 
imply $l(B)=1$, as the
following example shows.

Let $G$ be a finite group with a normal Sylow 
$p$-subgroup $P$.
Furthermore let $C_G(P) \leq P$. Thus $G$ is a $p$-solvable group 
with exactly one $p$-block (see \cite{Feit}, Chap. V, 3.11).
Since $ P$ is normal in $ G$ it is contained in the kernel of any 
irreducible
representation in characteristic $p$ . Thus we
have $ p \nmid \varphi(1)$ for all
$ \varphi \in \Irr_p(G)$.
A result of Fong (see \cite{Feit}, Chap. X, 3.2) yields 
$ \Phi(1) = |P|\varphi(1)$, where 
$\Phi$ denotes
the projective character corresponding to $\varphi$. Thus equality 
holds in (\ref{nr1e}), 
but $l(B)>1$ if $P \neq G$. 
\medskip

\noindent
{\bf 1.4 Blocks of $p$-solvable groups.}  
Let $B$ be a block of a $p$-solvable group 
$G$. We are going to show that the inequality (\ref{nr1e}) 
holds for $B$, which implies the Local Conjecture for blocks
of $p$-solvable groups. 

Let $\Phi_\varphi $ denote the projective
character corresponding to the irreducible Brauer character $\varphi$. 
Again by Fong's result we have  
$\Phi_\varphi(1) = |G|_p \varphi(1)_{p'}$
(see \cite{Feit}, Chap. X, 3.2) and, by a result of Brauer, we have 
$ \frac{|G|_p}{|D|} \mid \varphi(1)$ (see
\cite{Feit}, Chap. IV, 4.5). 
 Thus we get
$$ \mbox{{\Large $\frac{\dim B}{|D|}$}} = \mbox{{\Large 
$\frac{1}{|D|}$}} 
\sum_{\varphi \in \Irr_p(B)} 
\varphi(1) \Phi_\varphi(1) 
    =  \mbox{{\Large $\frac{|G|_p}{|D|}$}}  
\sum_{\varphi \in \Irr_p(B)} \varphi(1) 
\varphi(1)_{p'}
    \leq   \sum_{\varphi \in \Irr_p(B)} \varphi(1)^2.$$
For proving Theorem 1(b) it remains to characterize equality in (\ref{strong}).
From the discussions above we have
$$\mbox{{\Large $\frac{\dim B}{|D|}$}} 
\leq \sum_{\varphi \in \Irr_p(B)} \varphi(1)^2\le 
l(B)\max_{\varphi \in \Irr_p(B)} \{\varphi(1)^2\}~,$$
where the first inequality is an equality if and only if 
$\frac{|G|_p}{|D|}=\varphi(1)_p$ for all $\varphi\in\Irr_p(B)$.
Clearly, the second inequality becomes an equality if and only if
all irreducible Brauer characters have the same degree. This proves
Theorem 1(b). 
\medskip

\medskip

\noindent
{\bf 1.5 Connections to classical open problems.}
Let $k(B)$ denote the number of ordinary irreducible characters in
an arbitrary block  
$B$, with defect group $D$. 
Clearly $l(B) \leq k(B)$. Problem 20 of Brauer's famous list
\cite{Brauer} asks whether always $k(B)\le |D|$. 
If we assume that this has an affirmative answer we get $ l(B)|D| 
\leq |D|^2$, and our Local Conjecture
leads to 
\begin{eqnarray} \label{nr1c}
 \frac{\dim B}{|D|^2} \leq \sum_{\varphi \in \Irr_p(B)} \varphi(1)^2.
\end{eqnarray} 
This would obviously hold true if 
\begin{eqnarray} \label{nr1d}
\Phi_\varphi(1) \leq |D|^2 \varphi(1)  \  \mbox{ for all} \  \varphi \in
 \Irr_p(B)~,
\end{eqnarray}
where $\Phi_{\varphi}$ 
is the projective character corresponding to $\varphi$.
Unfortunately, there are examples for which (\ref{nr1d}) is wrong. 
For instance,
in characteristic $3$ we have $\Phi_\varphi(1) =99$ for the trivial character
$\varphi$ of the alternating group A$_7$ and  $|D|^2\varphi(1) =81$.
\medskip

\noindent
{\bf 1.6}
What are the connections between the Global Conjecture from
\cite{Willems} and the Local Conjecture?
An affirmative answer to the Local Conjecture would also prove a weaker
version of the Global Conjecture. In fact, if $l(G)$ denotes the number
of $p'$-conjugacy classes of $G$, then
by summing up (\ref{nr1b}) over all blocks $B$ of $G$ we would get
\begin{eqnarray} \label{nr1f}
\frac{|G|_{p'}}{l(G)} \leq  \sum_{\varphi \in \Irr_p(G)} \varphi(1)^2
\end{eqnarray}  
with equality if and only if $G$ is a $p$-group.
The latter fact can be seen as follows. The Local Conjecture implies
$$ |G| = \sum_B \dim B \leq \sum_B l(B)|D|  
\sum_{\varphi \in \Irr_p(B)} \varphi(1)^2 \leq
      l(G)|G|_p  \sum_{\varphi \in \Irr_p(G)} \varphi(1)^2.$$
Thus equality in (\ref{nr1f}) forces  equality in (\ref{nr1b}) 
for each block $B$, hence $l(B)=1$, and $l(B)=l(G)$. 
Thus the trivial character is the only irreducible
Brauer character and therefore $G$ has only one conjugacy class of
$p'$-elements.
This shows that $G$ is a $p$-group.
On the other hand, if $G$ is 
a $p$-group, then obviously equality holds in (\ref{nr1f}). 

(Note that in general $\Phi_\varphi(1) \not\leq
|P|l(G)\varphi(1)$, where $P$ is a Sylow $p$-subgroup of $G$. 
Take the projective character corresponding to the 
trivial character
for the smallest Suzuki group in characteristic $2$.)

\medskip

\noindent
{\bf 1.7} As mentioned in the introduction we do not know any example in
odd characteristic of a block $B$ for which
\begin{eqnarray} \label{slc}
 \frac{\dim B}{l(B)|D|} \leq \max_{\varphi \in \Irr_p(B)} \varphi(1)^2 
\end{eqnarray} 
 does not hold. So let us assume that $p \not=2$ and that (\ref{slc}) 
is true for
 all blocks. If we now sum up (\ref{slc}) over all blocks we get the 
Weak Global Conjecture
 \begin{eqnarray} \label{wgc}
  \frac{|G|_{p'}}{l(G)} \leq  \max_{\varphi \in \Irr_p(B)} \varphi(1)^2.
 \end{eqnarray}
 Moreover, and here $p\not=2$ is crucial for what follows, equality holds
 if and only if $G=O_p(G) \times O_{p'}(G) $, where $O_{p'}(G)$ is 
abelian. Thus the assumption (\ref{slc})
 for all blocks $B$
 leads to a complete proof of the Weak Global Conjecture in the case 
of odd characteristic.
 
 Indeed if $G=O_p(G) \times O_{p'}(G) $ where $O_{p'}(G)$ is abelian,
then equality holds in (\ref{wgc}).
In order to prove the converse suppose that 
equality holds in (\ref{wgc}). 
First we prove by induction that $G$ is $p$-solvable.
We put $H=O_{p'}(G)$ and $ \overline{G} =G/H$.
Since the group algebra of $\overline{G}$ is a direct 
summand of the group algebra of $G$
over the underlying field, one can show that
$$ \frac{|\overline{G}|_{p'}}{l(\overline{G})} = 
\max_{\varphi\in\Irr_p(\overline{G})} 
\{\varphi(1)^2\}.$$
If $H \not= \langle 1 \rangle$, then by induction 
$\overline{G}$ is $p$-solvable. Hence, $G$ is $p$-solvable.
Next we assume that $O_p(G) \not= \langle 1 \rangle$.
If $ \hat{G} = G/O_p(G)$, then equality in (\ref{wgc}) 
implies that
$$ \frac{|\hat{G}|_{p'}}{l(\hat{G})} = \max_{\varphi\in\Irr_p(\hat{G})} 
\{\varphi(1)^2\}~,$$
and the induction argument yields again that $\hat{G}$ is $p$-solvable.
Therefore, $G$ is also $p$-solvable.
Finally, we may assume that $G$ has a minimal normal subgroup, say $N$,
which is a direct product of isomorphic non-abelian simple groups whose
orders are divisible by $p$.
The equality in (\ref{wgc}) implies that 
all blocks of $G$ are of maximal defect. 
Thus, by Kn\"orr's theorem (see \cite{Navarro}, 9.26),
all blocks of $N$ are of maximal defect.  By a result of Michler 
(see \cite{Michler}, Theorem 5.3),
a finite non-abelian simple group has blocks of different 
defects for $p \not=2$.
Thus  we have a contradiction in this case, and the claim that $G$ 
is $p$-solvable
has been proved.

Next we show that $G$ is an abelian $p'$-group.
Since
$$ |G|_{p'} = \sum_{\varphi \in \Irr_p(G)} 
\varphi(1)_{p'}\varphi(1) \leq
l(G) \max_{\varphi \in \Irr_p(G)} \{ \varphi(1)^2\} =  |G|_{p'}~, $$
all $ \varphi$ are of degree $1$.
In particular
$$  G' \leq \cap_{ \varphi \in \Irr_p(G)} \kernel \varphi = O_p(G)~,$$
where $G'$ denotes the commutator subgroup of $G$. 
Thus $G$ has a normal Sylow $p$-subgroup
with abelian factor group. 
The equation
$$ \frac{|G|_{p'}}{l(G)} = \max_{\varphi\in\Irr_p(G)} 
\{\varphi(1)^2\} = 1 $$
shows that all $p'$-conjugacy classes have length $1$. 
Thus a Hall $p'$-subgroup of $G$
is central in $G$, which proves the assertion.


\section{Eigenvalues and Brauer characters} \label{Sec-eigen}


Associated to 
any $p$-block $B$ is the Cartan matrix
$C = C_B$.  If $ \langle \cdot \, , \cdot \rangle$ denotes 
the Euclidean inner product on $\R^l$ and
$ \overline{\varphi} = (\varphi_1(1), \ldots, \varphi_l(1))$
the dimension vector, then the dimension of the block  can be expressed as
$\dim B \, = \, \langle C \overline{\varphi}, 
\overline{\varphi}\rangle.$
If $\rho(C)$ denotes the Perron-Frobenius eigenvalue of $C$, i.e.
$$ \rho(C) = \max \, \{ |\lambda| \mid \lambda \  \mbox{is an eigenvalue 
of $C$} \},$$
then, since $C$ is a Hermitian matrix,
$$
  \frac{\langle C \overline{\varphi}, \overline{\varphi} 
\rangle}{\langle  \overline{\varphi}, \overline{\varphi} \rangle}
 \leq \rho(C).
$$ 
In particular we have the following result.

\begin{Proposition}  \label{Wada} {\rm (Kiyota-Wada \cite{Kiyota-Wada})}
\begin{equation} \label{nr2a}
  \frac{\dim \, B}{\rho(C)} = \frac{\langle C \overline{\varphi}, 
\overline{\varphi}\rangle}{\rho(C)} \ 
 \leq \  \langle  \overline{\varphi}, \overline{\varphi} \rangle \  = \  
\sum_{\varphi \in \Irr_p(B)} \varphi(1)^2.
\end{equation} 
\end{Proposition} 

In order to prove our Local Conjecture it therefore 
suffices to show that
$ \rho(C) \leq l(B)|D|$ holds for any block. 
Usually, $\rho(C)$ is hard to determine, and there is no general upper 
bound known which only depends on invariants of $B$, and not on the 
group. 
Moreover, equality in (\ref{nr2a}) does not seem to be easy to 
characterize. In fact, 
suppose equality holds in (\ref{nr2a}).
Since $C$ is Hermitian, this implies $ C \overline{\varphi} = 
\rho(C) \overline{\varphi}$.
Thus if $\Phi_\varphi$ denotes the projective character corresponding 
to $\varphi$,
then we have
$ \Phi_\varphi(1) = \rho(C) \varphi(1).$
By Theorem 3 of \cite{Brauer76} we have
$\gcd(\Phi_\varphi(1) \mid \varphi \in \Irr_p(B)) = p^au_B$,
where $|G|_p = p^a$ and 
$u_B$ is an integer not divisible by $p$. 
Moreover, by Theorem 2 of the same paper we have
$\gcd(\varphi(1) \mid \varphi \in \Irr_p(B)) = p^{a-d}u_B$,
where $|D|=p^d$.
Since  $\Phi_\varphi(1) = \rho(C) \varphi(1)$ we conclude that
$\rho(C) = |D|$.
On the other
hand if $ \Phi_\varphi(1) = |D| \varphi(1)$ for all $ \varphi \in 
\Irr_p(B)$, then
equality holds in (\ref{nr2a}).
\smallskip

Such situations occur for all blocks if $G$ is $p$-solvable of 
$p$-length $1$.
In particular, equality in (\ref{nr2a}) does not imply $l(B)=1$ 
in general.

Despite the difficulties in general, for blocks of $p$-solvable groups
we have the following characterization
of equality in (\ref{nr2a}). 

\begin{Proposition} Let $G$ be a $p$-solvable group. For the principal 
$p$-block $B$ the following statements are equivalent:
    \begin{itemize}
     \item[{(a)}] $\frac{\dim \, B}{\rho(C)}  = \  
             \sum_{\varphi \in \Irr_p(B)} \varphi(1)^2,$ i.e. 
          equality holds in {\em (}\ref{nr2a}{\em )}.
     \item[{(b)}] $G$ is of $p$-length $1$. 
      \end{itemize} 
\end{Proposition}

\noindent
{\it Proof.} As mentioned above, (b) implies (a) 
(see \cite{Schwarz}, Satz 4.2). 
To prove the other direction let $P$ denote a Sylow $p$-subgroup of $G$.
By the remarks preceding the proposition we have 
$\Phi_\varphi(1) = |P| \varphi(1)$
for all $ \varphi \in \Irr_p(B)$. This implies $ \Phi_\varphi = 
\varphi \Phi_1$ where
$\Phi_1$ corresponds to the trivial character. 
If $ p \mid \varphi(1)$, then
the multiplicity of the trivial character in 
$ \varphi \overline{\varphi}$
is at least $2$ (see \cite{Huppert}, Chap. VII, 8.5 d)). Thus
$$ 2 \leq (\Phi_1, \varphi \overline{\varphi}) = 
(\Phi_1\varphi, \varphi) =
            ( \Phi_\varphi, \varphi) = 1, $$
            a contradiction.
So all irreducible Brauer characters in $B$ are of $p'$-degree. 
Since $G$ is
assumed to be $p$-solvable and $B$ is the principal block,
we conclude that
$B$ is isomorphic to the group algebra of $G/O_{p'}(G)$. The condition 
$p \nmid \varphi(1)$
for all $ \varphi$ in $B$ implies that $G/O_{p'}(G)$ has a normal 
Sylow $p$-subgroup (a well-known
result of Wallace, see \cite{Wallace} or \cite{Brockhaus}). 
Thus $G$ is of $p$-length $1$. 
\smallskip

Our next aim is to establish a linear algebra criterion for proving
the Local Conjecture, namely Proposition 2 as stated in the Introduction.
We begin by recalling some useful fundamental facts about Cartan matrices
of blocks. 
            
\begin{Lemma} \label{pos} The Cartan matrix of any 
$p$-block $B$ is 
positive definite.
\end{Lemma} 

\noindent
{\it Proof.} 
If $D$ denotes the decomposition matrix of $B$, then $C=D^TD$. 
Then we have
$(Cx,x) = (D^TDx,x) = (Dx,Dx)$
for all $ x \in \R^l$, where $l=l(B)$. If $(Dx,Dx)=0$ then 
$Dx=0$ and therefore 
$Cx=0$. But $C$ is regular, hence $x=0$.
This proves that $(Cx,x)>0$ for all $ 0 \not= x \in \R^l$.

\begin{Corollary} \label{spur} If $C$ is the Cartan matrix of a 
$p$-block, then
               $$ \rho(C) \leq \tr C.$$
\end{Corollary}               

\noindent
{\it Proof.} Since $C$ is positive definite by Lemma \ref{pos},
all eigenvalues of $C$ are
positive, and they sum up to the trace of $C$.
In particular $ \rho(C) \leq \tr C$.

\begin{Proposition} \label{Prop-trace}
Let $B$ be a $p$-block of a finite group, with Cartan matrix $C$. Then
     $$ \frac{\dim B}{\tr C} \leq \sum_{\varphi \in \Irr_p(B)} 
\varphi(1)^2, $$
     and equality holds if and only if $l(B)=1$.
\end{Proposition}

\noindent
{\it Proof.} 
By Proposition \ref{Wada} and Corollary \ref{spur} we have
$$  \frac{\dim B}{\tr C} \leq  \frac{\dim B}{\rho(C)} \leq 
\sum_{\varphi \in \Irr_p(B)} \varphi(1)^2.$$
Suppose equality holds, i.e. $ \tr C = \rho(C)$. Since all eigenvalues
are positive, and their sum is $\tr C$, we deduce that $l(B)=1$. 

\smallskip

The above proposition gives a possible linear algebra approach to
proving the Local Conjecture, provided the following question has an
affirmative answer. 

\begin{Question} Let $B$ be a $p$-block of a finite group, 
with defect group $D$ and
Cartan matrix $C$. Is it true that 
  \begin{eqnarray} \label{nr2c}  \tr C \leq l(B) |D|? \end{eqnarray}
\end{Question}

\medskip

Since $C$ is positive definite we have
$$ \det C \leq \prod_{\varphi \in \Irr_p(B)} c_{\varphi \varphi}, $$
hence
$$ l \,  \sqrt[l]{ \det C} \leq l \,  \sqrt[l]{\prod_{\varphi \in 
\Irr_p(B)} c_{\varphi \varphi}} \leq
      \sum_{\varphi \in \Irr_p(B)} c_{\varphi \varphi}  = \tr C.
$$      
Thus the truth of (\ref{nr2c}) would imply 
$ \sqrt[l]{\det C} \leq |D|$. This is indeed the case, as one can see
using well-known properties of elementary divisors of $C$.


\section{Blocks with cyclic defect groups} \label{Sec-cyclic}


This section deals with blocks with cyclic defect group and contains 
a proof  of
Theorem 1(c).
The proof is very much 
based on the well-known structure theory of these blocks, and it
also follows from work of Kiyota and Wada \cite{Kiyota-Wada}.
For the sake of completeness and for the convenience of the reader 
we include a direct proof here.
We actually show the following stronger statement. 

\begin{Proposition} For a $p$-block 
$B$ with cyclic defect group $D$ we have
\begin{eqnarray} \label{cyclic}
   \frac{\dim \, B}{|D|}  
           \leq  \sum_{\varphi \in \Irr_p(B)} \varphi(1)^2 
\end{eqnarray}           
with equality if and only if the corresponding Brauer tree of $B$ 
is a star with
exceptional vertex, if there is one, in the center.
\end{Proposition}

\noindent
{\it Proof.} 
Let $e$ denote the inertial index of $B$. In other words $e$ is
the number $l(B)$ of irreducible Brauer
characters in $B$. If we put $ m = \frac{|D|-1}{e}$,
then the shape of the Brauer tree yields
$$  0 \leq C = (c_{i,j}) \leq  \left( \begin{array}{ccc}
         m+1 & & \\[2ex]
         m  & \ddots & m \\[2ex]
         &   & m+1    \end{array} \right) $$
for the Cartan matrix $C$ of $B$ where $ \leq $ for matrices means 
$ \leq$ in
each position. If we put $ \overline{\varphi} = 
(\varphi_1(1),\ldots, \varphi_e(1))$
and $ \overline{1} = (1,\ldots,1)$, then by the 
Cauchy-Schwarz inequality
$$ |\langle \overline{\varphi}, \overline{1} \rangle |^2 
\leq \ \langle \overline{\varphi}, \overline{\varphi}\rangle \ 
        \langle \overline{1}, \overline{1}\rangle,$$
hence $$ {\rm (*)} \qquad ( \sum_{i=1}^e \varphi_i(1))^2 \leq e 
\sum_{i=1}^e \varphi(1)^2.$$
Thus we get
$$ \begin{array}{rcl}
    \dim B & = & \sum_{i=1}^e \varphi_i(1) \Phi_i(1) \\[1ex]
           & = &  \sum_{i,j=1}^e \varphi_i(1) c_{i,j} \varphi_j(1) \\[1ex]
           & \leq & m \sum_{i,j=1}^e \varphi_i(1) \varphi_j(1) +
                \sum_{i=1}^e \varphi_i(1)^2  \\[1ex]
           & = & m \sum_{i=1}^e \varphi_i(1) \sum_{j=1}^e \varphi_j(1) +
                \sum_{i=1}^e \varphi_i(1)^2  \\[1ex]              
           & \leq & me \sum_{i=1}^e \varphi_i(1)^2 + 
                       \sum_{i=1}^e \varphi_i(1)^2 \qquad 
\qquad \qquad \mbox{by} \  (*)  \\[1ex]
           & = & (me+1)\sum_{i=1}^e \varphi_i(1)^2 \\[1ex]
           & = & |D| \sum_{i=1}^e \varphi_i(1)^2.
   \end{array}
$$
Suppose equality holds in (\ref{cyclic}). It follows that
$$   C = (c_{i,j}) =  \left( \begin{array}{ccc}
            m+1 & & \\
             m  & \ddots & m \\[2ex]
                &     & m+1    \end{array} \right). $$
This means that the Brauer tree is a star with the exceptional 
vertex (if there is one)
in the center. 


\section{ Cartan matrices of tame blocks} \label{Sec-tame}


In this section we give a proof of our Local Conjecture for blocks
of tame representation type, i.e. for all $2$-blocks of finite 
groups having a dihedral, semidihedral or generalized quaternion defect 
group. For the proof we use the linear algebra approach outlined
in Section \ref{Sec-eigen}, namely we are going to show that 
the trace of the Cartan matrix can be bounded above in terms 
of the order of the defect group. The
Cartan matrices of tame blocks are known due to the extraordinary 
classification of K. Erdmann of tame blocks up to Morita equivalence
\cite{Erdmann}. 

Blocks of tame representation type have at most three simple modules. 
Since our conjectures are true 
for blocks with one simple module,
we can restrict our attention to the blocks with two or
three simple modules.  
We refer to \cite{Erdmann} for background on tame blocks,
and also for details and notations. 

\begin{Proposition} \label{tame1} Let $B$ be a block of a finite group
of tame representation type, with defect group $D$. 
Then for the trace of the 
Cartan matrix $C$ of $B$ the following holds:
\begin{enumerate}
\item[{(a)}] If $B$ has two simple modules, then we 
have $\ \tr C \, \le \frac{5}{4}|D|+2.$ 
\item[{(b)}] If $B$ has three simple modules, then we have 
$\tr C \, \le \frac{3}{2}|D|+4$.
\end{enumerate}
\end{Proposition}

\noindent
{\it Proof.} 
(a) From Erdmann's classification we see that such a block $B$ has a
Cartan matrix of the form 
$$\left( \begin{array}{cc} 4k & 2k \\ 2k & k+r 
\end{array} \right)
$$
with natural numbers $k$ and $r$, where 
$\{k,r\}=\{1,\frac{|D|}{4}\}$ or $\{k,r\}=\{2,\frac{|D|}{4}\}$.
For any possiblity of $k$ and $r$ we obtain that 
$\tr C \, \le \frac{5}{4}|D|+2$ (note that $|D|\ge 8$, since 
a block with Klein four defect group cannot have two simple modules).
\smallskip

(b) There are various
families of algebras in Erdmann's classification
which (may) contain blocks of finite groups. They are denoted 
$D(3{\mathcal A})_1$, $D(3{\mathcal B})_1$ and $D(3{\mathcal K})$
in the dihedral case, $SD(3{\mathcal A})_1$, $SD(3{\mathcal B})_1$,
$SD(3{\mathcal C})_2$ and $SD(3{\mathcal D})$ in the semidihedral case,
and $Q(3{\mathcal A})_2$, $Q(3{\mathcal B})$ and $Q(3{\mathcal K})$
in the quaternion case. 

The families $D(3{\mathcal A})_1$ and $Q(3{\mathcal A})_2$ have Cartan 
matrices of the form 
$$\left( \begin{array}{ccc}
4k & 2k & 2k \\ 2k & k+a & k \\ 2k & k & k+a
\end{array} \right)~,
$$
where $k\in\mathbb{N}$ and $a\in\{1,2\}$. 
Blocks of finite groups only occur for $k=\frac{|D|}{4}$. 
 In particular, $\tr C \, \le \frac{3}{2}|D|+4$. 
\smallskip

The algebras $D(3{\mathcal K})$ and $Q(3{\mathcal K})$ have Cartan 
matrices of the form 
$$\left( \begin{array}{ccc}
2a & a & a \\ a & k+a & k \\ a & k & k+a
\end{array} \right)~,
$$
where $k\in\mathbb{N}$ and $a\in\{1,2\}$. 
For blocks of finite groups we have $k=\frac{|D|}{4}$. In particular
$\tr C \, \le \frac{|D|}{2}+8\le \frac{3}{2}|D|+4$ (since $|D|\ge 4$
for tame blocks). 
\smallskip

The algebras in the families $SD(3{\mathcal B})_1$,
$SD(3{\mathcal D})$ and $Q(3{\mathcal B})$ all have Cartan matrices
of the form
$$\left( \begin{array}{ccc}
4 & 2 & 2 \\ 2 & s+1 & 1 \\ 2 & 1 & 3
\end{array} \right)
$$
for a natural number $s\ge 2$.  
Blocks of finite groups only occur 
for $s=\frac{|D|}{4}\ge 2$. In particular we get
that
$\tr C \, \le \frac{|D|}{4}+8\le \frac{3}{2}|D|+4.$ 
\smallskip

The algebras in the family $D(3{\mathcal B})_1$ have Cartan matrix
of the form
$$\left( \begin{array}{ccc}
4 & 2 & 2 \\ 2 & s+1 & 1 \\ 2 & 1 & 2
\end{array} \right)
$$
for some $s\ge 1$.  
Blocks of finite groups occur for the parameter
$s=\frac{|D|}{4}$, hence $\tr C \, \le \frac{|D|}{4}+7$ for blocks. 
\smallskip

The algebras in the family $SD(3{\mathcal A})_1$ have Cartan matrix
of the form
$$\left( \begin{array}{ccc}
4k & 2k & 2k \\ 2k & k+1 & k \\ 2k & k & k+2
\end{array} \right)
$$
for some natural number $k$.  Blocks of finite groups occur 
for parameters $k=\frac{|D|}{4}$, i.e. $\tr C \, \le \frac{3}{2}|D|+3$
for blocks in this family.
\smallskip

There are two more families for which it is still an open question 
whether blocks of finite groups actually occur in these families. 

The Cartan matrices of algebras in the
family $SD(3{\mathcal H})$
are of the form
$$\left( \begin{array}{ccc}
3 & 2 & 1 \\ 2 & s+2 & s \\ 1 & s & s+1
\end{array} \right)
$$
for some $s\ge 2$. If blocks occur, then they occur for 
$s=\frac{|D|}{4}$. 
In particular, for blocks we have
$\tr C\le \frac{|D|}{2}+5\le \frac{3}{2}|D|+4.$ 
\smallskip

The Cartan matrices for algebras in the family $SD(3{\mathcal C})$
are of the form
$$\left( \begin{array}{ccc}
k+s & k & k \\ k & k+1 & k-1 \\ k & k-1 & k+1
\end{array} \right)~.
$$
If blocks of finite groups occur, then we have
$\{k,s\}=\{2,\frac{|D|}{4}\}.$ For $k=2$, $s=\frac{|D|}{4}$ 
we get $\tr C\le \frac{|D|}{4}+8$, in the other case we get
$\tr C \, \le \frac{3}{4}|D|+4.$ 

This completes the proof of the proposition. 
\medskip

As an immediate consequence of the above proposition,
we obtain that if $B$ is a block of a finite group
with tame representation type, then $\tr C\le l(B)|D|$, where 
$D$ is the defect group. 
(Note that a block with a Klein four group as a defect group 
has three irreducible Brauer characters.) By Proposition
\ref{Prop-trace} we can conclude that the Weak Conjecture holds
for all tame blocks, thus proving Theorem 3.


\end{document}